\documentclass{article}
\def\uu{\bigsqcup}

\def \Z {{\mathbf {Z}}}

\title{ Жесткие  пуассоновские надстройки без корней}
\author{ В.В. Рыжиков}
\date{}

\usepackage[T1]{fontenc} 
\usepackage[cp1251]{inputenc}  
\usepackage[russian]{babel}
\begin{document}

\maketitle

\begin{abstract}   

Examples of \bf rigid Poisson suspensions without roots \rm are presented. The discrete rational component in spectrum of an ergodic automorphism S prevents some roots from existing. If S is tensorly multiplied by an ergodic automorphism of the space with a sigma-finite measure, discrete spectrum disappears in this product, but like the smile of Cheshire Сat, the memory of it can remain in the form of the absence of roots. In additional conditions, this effect is inherited by the Poisson suspension over the above product.

\vspace{2mm}
Предъявлены  примеры жестких   пуассоновские надстройки  без корней. Наличие в спектре  эргодического автоморфизма $S$ дискретной  рациональной компоненты несовместимо с наличием определенных корней у него.  Если $S$ тензорно умножить на эргодических автоморфизм  пространства с сигма-конечной мерой, в произведении этот дискретный спектр исчезает,  но подобно улыбке Чеширского Кота    остатется память о нем   в виде  отсутствия корня. При некоторых дополнительных условиях этот эффект наследует пуассоновская надстройка над  произведением.

\vspace{2mm}
\it Ключевые слова:  \rm пуассоновские  надстройки, гауссовские автоморфизмы,  жесткость,  корни автоморфизма,  спектр.

\end{abstract}
 
\Large
\section{Введение}
Пусть $T$ -- автоморфизм пространства $(X,\mu)$ с сигма-конечной мерой
$\mu$, 
отвечающий ему ортогональный оператор в $L_2(\mu)$ мы также обозначаем  через $T$.
Группа всех автоморфизмов пространства $(X,\mu)$ вкладывается  непрерывно в группу автоморфизмов  пространства конфигураций с вероятностной пуассоновской мерой (см. \cite{Ne},\cite{Ro}). Образ автоморфизма $T$ пуассоновском вложении при называется пуассоновской надстройкой и обозначается через $P(T)$. 
Надстройка $P(T)$ над автоморфизмом $T$ обладает спектральным двойником:  гауссовским автоморфизмом  $G(T)$, ассоциированным  с действием ортогонального оператора $T$ на пространстве с вероятностной гауссовской мерой.  Автоморфизмы $G(T)$ и $P(T)$ спектрально изоморфны, но могут иметь существенно  различные метрические (от слова мера) свойства, включая алгебраические свойства. 

Хорошо известно, что надстройки   $P(T)$ и  $G(T)$  как операторы изоморфны оператору 
$$1\oplus \bigoplus_{n=1}^\infty  T^{\odot n},$$ 
$1$ -- тождественный оператор в одномерном пространсте, $T^{\odot n}$   -- симметрическая $n$-степень оператора $T$.
Гауссовские автоморфизмы обладают корнями, их центализатор  обширен, так как содержит инъективный образ  обширного  централизатора  оператора $T$ в группе  всех ортогональных операторов на вещественном пространстве $L_2(\mu)$. 

Пусть $(\perp)$ обозначает следующее свойство:    \it спектральная мера автоморфизма    взаимно сингулярна с ее сверточными степенями.   \rm Это эквивалентно тому, что оператор $T$ не имеет ненулевого сплетения с тензорными степенями $T^{\otimes n}$ при $n>1$. Если централизатор автоморфизма $T$ в группе всех автоморфизмов пространства $(X,\mu)$  состоит только  из степеней автоморфизма  и $T$ обладает свойством $(\perp)$, то пуассоновская    надстройка $P(T)$ также обладает тривиальным централизатором. Этот факт является   следствием    предложения 5.2 работы Э. Руа \cite{Ro} о том, что  $ C(P(T))= P(C(T))$ -- централизатор  надстройки $P(T)$ есть надстройка над централизатором  автоморфизма $T$.

Если у $T$ нет корней и $T$ обладает свойством $(\perp)$, то $P(T)$ также не имеет корней, это контрастирует с  тем, что   всякий эргодический гауссовский  автоморфизм $G(T)$ входит  в  континуум неизоморфных потоков, следовательно,  имеет континуум корней.  Другие впечатляющие контрасты см. в \cite{PR}. Отметим также, что примеры  перемешивающих надстроек с тривиальным централизатором  обеспечивают работы  \cite{DR},\cite{RT},\cite{R24}. Напомним, что $P(T)$ обладает перемешиванием, если    $T^n\to 0$ при $n\to\infty$.  В  \cite{R24}  предъявлены  автоморфизмы с сингулярным спектром, сверточные степени которого лебеговские, что влечет за собой свойство  $(\perp)$.
 Неперемешивающие  пуассоновские надстройки $P(T)$ с тривиальным централизатором  получаются из конструкций нежестких автоморфизмов, предложенных в  \cite{R19}.  

Если автоморфизм $T$  жесткий, т.е. $T^{n_i}\to I$  для некоторой
последовательности $n_i\to\infty$,    централизатор $C(P(T))$  континуален. Но если автоморфизм   $T$ не имеет  корней, но обладает свойством $(\perp)$, надстройка  $P(T)$ тоже не имеет корней. 
Хорошо известно, что типичный автоморфизм вероятностного пространства слабо перемешивает (это свойство означает  непрерывность спектра  автоморфизма при его действия на  на пространство, ортогональное к констатнтам) и обладает жесткостью.  Дж. Кинг доказал, что   типичный автоморфизм обладает корнями  всех степеней \cite{K}. 

 \it  Задача построения жесткого слабо перемешивающего автоморфизма без корней \rm требует значительной изобретательности.  На возможность таких построений  методом косых произведений  намекала давняя работа  А.М. Степина \cite{St} об   автоморфизмах без квадратных корней. 
Т. Адамс    сообщил автору, что  подходщим  примером является жесткий автоморфизм ранга один, предложенный А. дель Джунко и Д. Рудольфом в \cite{JR}. Простая модификация этого примера  обеспечивает бесконечность меры фазового пространства, сохраняет свойства жесткости и $(\perp)$, но доказательство   отсутствия корней, вероятно,  требует  привлечения тонкого анализа.   Мы предлагаем относительно простое  решение упомянутой задачи  в классе  пуассоновских надстроек.  Поясним
суть нашего подхода.

\bf Рациональный спектр. \rm Если эргодический автоморфизм $S$ вероятностного пространства $(X,\mu)$ как оператор обладает собственным числом $-1$,  то $S$ не  имеет  квадратного корня. Покажем это. Модуль собственной функции $f$ инвариантен относительно $T$,  в силу эргодичности $S$   постоянен,  пусть он равен $1$. Определим  
$$Y=\{x\in X: f(x)= e^{it}, 0\leq t <\pi\}.$$
Можно считать, что  $f=2\chi_Y-1$. Так как  $Rg$ --  тоже собственная  функция для $S$, а собственное значение $-1$  однократно,   получим
 $Rf=\pm f$. Тогда $Sf=R^2f=f$. Противоречие.

\bf Эффект Чеширского Кота. \rm  Мы напомнили хорошо известный факт о том, что дискретный рациональный спектр эргодического автоморфизма является препятствием к извлечению соответствующих корней из автоморфизма.
Пусть $\tilde T=T\otimes S$,  где $T$ --  автоморфизм пространства с сигма-конечной мерой, причем все степени $T^p$, $p>1$, эргодические,   а $S$   обладает рациональным  спектром. Тогда у произведения $\tilde T$  рациональный  спектр подобно Чеширскому Коту пропадает, но может остаться его улыбка в виде  отсутствия корня у автоморфизма $\tilde T$. При условии $(\perp)$ этот эффект распространяется на пуассоновскую надстройку $P(\tilde T)$.

\section{Надстройки,  не обладающие  корнями степени $\bf p$}  Пусть $T$  -- эргодический автоморфизм стандартного пространства $(Y, m)$ с сигма-конечной мерой (в качестве   $Y$ берем прямую линию с мерой Лебега $m$). Рассмотрим  отображение 
$$F(T,p): (Y\times \Z_p) \to  (Y\times \Z_p),$$ где  $\Z_p$ --  аддитивная группа вычетов мощности $p>1$,      образом: 
$$  F(T,p)(x,z):= (Tx, z+1), \ \ z=0,$$ 
$$   F(T,p)(x,z):= (x, z+1), \ \ \ z\neq 0.$$ 
\\ Положим $\mu=m\times \nu$, где $\nu$ -- равномерная вероятностная мера (Хаара) на $\Z_p$. Преобразование $F(T,p)$, очевидно,  сохраняет меру $\mu$. 

\vspace{2mm}
\bf Лемма 1. \it Если $T$  -- эргодический автоморфизм, то  преобразование $F_p=F(T,p)$ не имеет корней степени $p$. \rm

\vspace{2mm} Доказательство. 
Пусть $R^p=F_p$.  Степень  $F_p^p$  имеет  $p$  эргодических компонент, рассмотрим  $X=\uu_0^{p-1} Y_i$ --  разбиение всего пространства на соответствующие инвариантные множества, на которых степень эргодична. Автоморфизм  $R$ коммутирует с $F_p^p$,  поэтому множества $RY_i$ также инвариантны относительно $F_p^p$  и на этх множествах  степень $F_p^p$  действует эргодично. Но отсюда вытекает, что для некоторого $k$ выполнено
$$RY_0=Y_k=F_p^kY_0, \ \ R^pY_0=F_p^{kp}Y_0=Y_0,$$ 
что влечет за собой $ R^p\neq F_p.$  Лемма доказана.

\vspace{2mm}
А.Б. Катком и А.М. Степиным  в \cite{KS} рассматривалась сильная сходимость степеней автоморфизма к оператору $-I$ на подпространстве в $L_2$, что влечет за собой  отсутствие подчинения четных сверток его спектру  и тем самым дает контпример к гипотезе 
А.Н. Колмогорова о  групповом свойстве спектра автоморфизмов пространства Лебега. 
 
 В.И. Оселедец в работе \cite{O} предложил использовать   неунитарные слабые пределы $aI$, $0<a<1$ (здесь предел рассматривается в пространстве, ортогональном константам). Наличие таких пределов  вынуждает    взаимную сингулярность всех сверточных степеней спектра автоморфизма,  что   обеспечивает выполнение свойства $(\perp)$. В  \cite{O} была высказана гипотеза о том, что соответствующие примеры можно найти в перекладываниях конечного числа отрезков, позднее автор  подтвердил  гипотезу Оселедца.

\vspace{2mm}
\bf Лемма 2. \it (i) Найдется автоморфизм $T:Y\to Y$ такой, 
  что для некоторого $a$, $0<a<1,$ и последовательностей $m_i, n_j\to\infty $   выполнено 
$T^{m_i}\to aI,$  $T^{n_j}\to I.$

(ii) Если    $T$ удовлетворяет (i), то автоморфизм $F_p=F(T,p)$ обладает свойством $(\perp)$.\rm

\vspace{2mm}
Доказательство.  Требуемые в (i) автоморфизмы $T$  легко строятся в классе преобразований ранга один (см. \S 4): на одной последовательности этапов 
обеспечивается  предел $aI$, а сходимость к $I$ реализуется на другой последовательности. 
Свойство (ii)  вытекает из сходимости 
$F_p^{pm_i}\to aI$ (здесь $I$ обозначает тождественный оператор
в $L_2(X,\mu)$). 
Поясним, как установить   свойство $(\perp)$.  Заметим, что  у операторов $F_p$ и $F_p\otimes F_p$ нет ненулевого сплетающего оператора. Действительно, если
$F_p J=J(F_p\otimes F_p),$ то  $$ F_p^{pm_i} J=J(F_p\otimes F_p)^{pm_i}, \ \ aJ=a^2J, \ \ J=0.$$
Аналогично устанавливается (при помощи  равенств  $ a^m J=a^nJ$)  дизъюнктность других сверточных степеней спектральной меры. Лемма доказана.

Таким образом, мы показали,  что   $F_p$, обладая  свойствами жесткости и  $(\perp)$,  не имет   корней степени $p$. С учетом    предложения 5.2  \cite{Ro}  приходим к следующему факту.

\vspace{2mm}
\bf  Теорема 1.    \it Для всякого $p>1$  жесткая пуассоновская надстройка $P(F_p)$ не имеет корня  степени $p$.\rm

 \section{Жесткие надстройки без  корней} 
\bf  Теорема 2.  \it Пусть $S$ --  автоморфизм с  рациональным  спектром, образующим группу $\{\exp(2k\pi i/p_1p_2\dots p_n)\}$, где $p_1, p_2,\dots,p_n$  -- всевозможные наборы различных простых чисел. Если для жесткого автоморфизма $T$ 
пространства с сигма-конечной мерой для всякого простого  $p>0$  степень $(S\times T)^p$ имеет в точности $p$ эргодических компонент и для некоторой последовательности $m_i$ выполнено  $T^{m_i}\to aI$, $0<a<1$, то пуасоновская надстройка $P(S\times T)$ является жесткой и не имеет корней.\rm 

\vspace{2mm}
Доказательство. 

\bf Отсутствие корней. \rm  Повторяя рассуждения, используемые в доказательстве теоремы 1, получаем, что автоморфизм $\tilde T= S\times T$ не имеет корней степени $p$ для всех простых чисел $p$, следовательно, не имеет никаких корней.

\bf Жесткость произведения  $ \bf S\times T$ . \rm  Автоморфизмы с одинаковым дискретным спектром изоморфны, автоморфизм  $S$ 
можно отождествить с  произведением $\prod_p  S_p$,  где $S_p$ -- сдвиг на $\Z_p$, а $p$ пробегает все простые числа.  Если $T^{h_i}\to I$ (здесь имеет место сильная операторная сходимость), то $T^{mh_i}\to I$,
поэтому для  всякого фиксированного $n$ имеем 
$$  T^{p_1p_2\dots p_n h_{i}}\ \to \ I,  \ i\to\infty.$$
Для автоморфизма $S$ для любой целочисленной последовательности $m_n$ выполнено
$$  S^{p_1p_2\dots p_n m_n}\ \to \ I,  при \ n\to\infty.$$
Cледовательно, для некоторой последовательности $i(n)\to\infty$ получаем сходимость
$ (S\times T)^{m_{i(n)}}\ \to I,$ которая влечет за собой  сходимость $ P(S\times T)^{m_{i(n)}}\ \to \ I$
(тождественные операторы, действующие   в разных пространствах, мы обозначили одинаково).

\bf Cвойство  $ \bf (\perp)$. \rm  
Из  $T^{m_i}\to aI$, $0<a<1$, вытекает, что  $z^{m_i}\to_w a$ в пространстве $L_2(\sigma)$, где $\sigma$ --  спектральная мера автоморфизма $T$ (подразумевается, что $\sigma$ -- мера на единичной окружности в комплексной плоскости, являющаяся мерой максимального спектрального типа для $T$).  Спектральная мера произведения $S\times T$ является сверткой спектральных мер сомножителей $S$ и $T$. Спектр $S$ дискретный, поэтому свертка
является взвешенной суммой сдвигов  меры $\sigma$. Запишем эту сумму в виде $\tilde \sigma=\sum_q c_q \sigma_q$.  Заметим, что $\sigma_q$ взаимно сингулярна с каждой  сверткой вида
$$ \tau=\sigma_{q_1}\ast \sigma_{q_2}\ast \dots \ast \sigma_{q_n}.$$ 
Это следует из того, что на $L_2(\sigma_q)$ для некоторой последовательности $i'\to\infty$
имеем слабую сходимость функций  $z^{m_{i'}}$ от комплексной переменной $z$, $|z|=1$, вида  $$z^{m_{i'}}\to_w a\lambda,$$ для некоторого $\lambda$,   $|\lambda|=1$,
причем в  пространстве  $L_2(\tau)$ наблюдаем сходимость 
 $$z^{m_{i'}}\to_w a^n\lambda_{1}\lambda_{2}\dots  \lambda_{n}, \
 |\lambda_{p}|=1,\ i'\to\infty. $$
  Очевидно, что меры $\sigma_{q}$ и $\tau$ не имеют общей части, иначе  $a=a^n, n>1$,  а это не так. 
Мы показали, что спектральная мера $\tilde \sigma$ сингулярна относительно своих сверточных степеней. А это означает выполнение свойства $(\perp)$ для автоморфизма $\tilde T$.

Таким образом, пуассоновская надстройка наследует свойство жесткости, отсутствие корней и является слабо  перемешивающей (имеет непрерывный спектр). Теорема доказана.

\section{Примеры подходящей базы для пуассоновской надстройки} 
Нам следует предъявить  автоморфизмы $T$, удовлетворяющие условиям теоремы 2. Напомним определение конструкции $T$ ранга один.

Фиксируем натуральное  число $h_1\geq 1$ (высота башни на этапе $j=1$), последовательность  $r_j\to\infty$ (число колонн, на которые виртуально разрезается башня этапа $j$)   и последовательность целочисленных векторов (параметров надстроек)   
$$ \bar s_j=(s_j(1), s_j(2),\dots, s_j(r_j-1),s_j(r_j)).$$  
На шаге $j=1$ задан  полуинтервал $E_1$. Пусть на  шаге $j$  определена  
 система   непересекающихся полуинтервалов 
$C_j, TC_j,\dots, T^{h_j-1}C_j,$
причем на $C_j,  \dots, T^{h_j-2}C_j$
пребразование $T$ является параллельным переносом. Такой набор   полуинтервалов  называется башней этапа $j$, их объединение обозначается через $X_j$ и тоже называется башней.

Представим   $C_j$ в виде   дизъюнктного объединения  полуинтервалов 
$C_j^i$, $i=1,2,\dots, {r_j},$ одинаковой длины.  
Для  каждого $i=1,2,\dots, r_j$ рассмотрим так называемую колонну  
$$C_j^i, TC_j^i ,T^2 C_j^i,\dots, T^{h_j-1}C_j^i.$$
К каждой  колонне с номером $i$  добавим  $s_j(i)$  непересекающихся полуинтервалов (этажей)  длины, равной длине интервала $C_j^i$.
Полученные  наборы интервалов  при  фиксированных $i$,$j$ называем надстроенными колоннами  $X_{i,j}$. Отметим, что при фиксированном $j$ по построению колонны $X_{i,j}$   не пересекаются. Используя параллельный перенос интервалов, преобразование $T$ теперь  доопределим так, чтобы колонны $X_{i,j}$ имели вид 
$$C_j^i, TC_j^i, \dots,  T^{h_j}C_j^i, T^{h_j+1}C_j^i, \dots, T^{h_j+s_j(i)-1}C_j^i,$$
  а   верхние этажи   колонн  $X_{i,j}$ ($i<r_j$) преобразование $T$ параллельным переносом  отображало в нижние
этажи колонн $X_{i+1,j}$: 
$$T^{h_j+s_j(i)}C_j^i = C_j^{i+1}, \ 0<i<r_j.$$ 
Положив $C_{j+1}= C^1_j$, замечаем, что все указанные этажи надстроенных колонн в новых обозначениях имеют вид 
$$C_{j+1}, TC_{j+1}, T^2 C_{j+1},\dots, T^{h_{j+1}-1}C_{j+1},$$
 образуя башню  этапа $j+1$ высоты  
 $$ h_{j+1} =h_jr_j +\sum_{i=1}^{r_j}s_j(i).$$

Частичное определение преобразования $T$ 
на этапе $j$ сохраняется на всех последующих этапах. 
В результате   получаем пространство  $X=\cup_j X_j$ 
и обратимое преобразование $T:X\to X$, сохраняющее  
стандартную меру Лебега $\mu$  на $X$.

Теперь зададим   параметры конструкции $T$, чтобы получить нужные свойства произведения $S\times T$.

\bf Произведение $\bf S^p\times T^p$ имеет  $\bf  p$ эргодических компонент.
\rm На  некотором множестве  этапов $J$ ( $J$ выбираем произольно с условием, что оно и его дополнение бесконечны)  мы обеспечим   свойство  жесткости конструкции $T$ и свойство произведения 
 $S^p\times T^p$ иметь ровно  $p$ эргодических компонент. Последнее  эквивалентно эргодичности  произведения $\hat S \times T^p $,
где $\hat S= \prod_{p'\neq p}S_{p'}$, а индекс $p'$  пробегает все простые числа, кроме $p$, так как  $S^p$ изоморфно прямой сумме $p$ копий, подобных  автоморфизму $\hat S$. 
Предположим, что  для всех   $j\in J$ числа $h_j$  взаимно просты с произведением $p_1\dots p_n$.  Это свойство легко обеспечить, так как на этапе $j-1$  мы можем определить $h_j$  любым числом, начиная с $2h_{j-1}$.

Положим $s_j(i)=0$ для всех $i=1,2,\dots, r_j$,  $r_j\to\infty$.    Тогда $$\mu(T^{h_j-n}C_j\Delta T^{p-n} C_j)/\mu(C_j)\to 0, \ n=1,2,\dots, p-1. \eqno (1)$$
Степень $T^p$  мало отличается от циклической перестановки этажей башни 
этого этапа, когда под действием $T$  верхний этаж целиком попадал бы в  нижний (этот случай имел бы место, если параметры $s_j(i)$ начиная с текущего  этапах были равны 0). 
Заметим, что   $\hat S$ является в точности циклической  престановкой из $p_1\dots p_n/p$  множеств $D_j, \hat SD_j,\dots, \hat S^{p_1\dots p_n/p-1}$ (считаем, что $n$ достаточно большое и $p$ находится среди $p_1,\dots,p_n$).
В силу того, что $h_j$ взаимно просто с $p_1\dots p_n/p$  и с учетом
$(1)$ получаем, что  $\hat S\times T^p$ мало отличается от циклической перестановки $h_jp_1\dots p_n/p$ множеств, являющихся декартовыми произведениями  $E_{q,h}=\hat S^qD_j\times T^hC_j$, $0\leq q< p_1\dots p_n/p$,
$0\leq h<h_j$.  Всякое множество (речь идет о подмножестве
пространства, на котором действует $\hat S\times T$), имеющее  конечную меру, аппроксимируется объединениями множеств $E_{q,h}$, $0\leq q< p_1\dots p_n/p$,
$0\leq h<h_j$. Автоморфизм с такой 
циклической  аппроксимацией   эргодический (аргументируем в духе работы \cite{KS}).

\bf Жесткость автоморфизма $\bf T$.  \rm Вытекает из того, что для   
  множеств  $A,B$, состоящих из этажей  башни этапа $j_0$,  для всех
этапов $ j>j_0$ для  $j\in J$ выполнено
$$|\mu(A\cap T^{h_j}B) -\mu(A\cap B)|\leq (\mu(A)+\mu(B))/r_j, $$
что влечет за собой  
$$T^{h_j}\to I, j\in J, \ j\to\infty.$$

\bf Свойство $\bf T^{m_i}\to aI$, $\bf 0<a<1$.  \rm    Ограничимся случам $a=1/2$. Для  $j\in J'$, (множество этапов $J'$ не пересекается с выбранным ранее множеством $J$) положим $r_j=2j$,   $s_j(i)=0$ при $i=1.2. \dots j$  и  $s_j(i)=h_j$ при  $i=j+1, j+2,\dots, 2j$.  
Тогда непосредственная проверка дает 
$$\mu(A\cap T^{h_j}B)= \mu(A\cap B)/2,\ j\in J',$$
следовательно,  
$$ T^{h_{2j}}\to_w I/2, \ j\to\infty,  j\in J'.$$

\bf Мера фазового пространства  $\bf X$  бесконечна. \rm  Очевидно, так как на  этапах $j\in J'$ к башне $X_j$ прибавлялось множество меры $\mu(X_j)/2$.

\newpage
Мы реализовали нужный нам эффект в классе жестких пуассоновских надстроек.   В заключение отметим, что теория пуассоновских надстроек
в свою очередь может предложить  задачу на экспорт.

\vspace{2mm}
\bf  Задача о слабо гомоклиническом элементе для жесткого слабо 
перемешивающего автоморфизма. \rm  
Эргодические пуассоновские надстройки $P$  (и  гауссовские автоморфизмы) обладают весьма специфическим свойством \cite{R}: найдется эргодический автоморфизм  $S$, для которого  
$$ \frac 1 N \sum_{n=1}^N  P^{-n}SP^n  \ \to \ I, \ n\to\infty.$$ 
 Обладают ли типичные автоморфизмы этим свойством, неизвестно.
Примеры  слабо перемешивающих жестких автоморфизмов
без этого свойства также неизвестны. Возможно, подходящим 
является преобразование  из \cite{JR}.

\vspace{2mm}
\bf Благодарности. \rm  Автор признателен Т. Адамсу, С. Еловацкому
 и \\ И. Подвигину за замечания и вопросы, послужившие поводом для написания этой заметки.

\end{document}